\documentclass{emsprocart}

\contact[geordie@mpim-bonn.mpg.de]{Geordie Williamson, Max-Planck-Institut f\"ur
  Mathematik, Vivatsgasse 7, 53111 Bonn, Germany}

\numberwithin{equation}{section}

\newtheorem{theorem}{Theorem}[section]

\newtheorem{conjecture}[theorem]{Conjecture}

\theoremstyle{definition}

\newtheorem{example}[theorem]{Example}
\newtheorem{remark}[theorem]{Remark}

\newcommand{\RM}{\mathbb{R}}
\newcommand{\QM}{\mathbb{Q}}
\newcommand{\ZM}{\mathbb{Z}}
\newcommand{\CM}{\mathbb{C}}

\newcommand{\reg}{\textrm{reg}}
\newcommand{\ample}{\textrm{ample}}
\newcommand{\simto}{\stackrel{\sim}{\to}}
\newcommand{\into}{\hookrightarrow}

\newcommand{\Hom}{\operatorname{Hom}}

\DeclareMathOperator{\End}{End}
\newcommand{\gr}{\operatorname{gr}}

\newcommand{\id}{\textrm{id}}
\newcommand{\hg}{\mathfrak{h}}
\newcommand{\bim}{R\textrm{-bim}}
\newcommand{\bimf}{R\textrm{-bim}^{fg}}

\newcommand{\HC}{\mathcal{H}}
\newcommand{\mini}{\textrm{min}}

\DeclareMathOperator{\tr}{tr}

\DeclareMathOperator{\ch}{ch}

\newcommand{\un}{\underline}

\newcommand{\hpos}{\hg_\reg^+}
\newcommand{\hspos}{\hg^{*+}_\reg}

\title[The Hodge theory of the Hecke category]{The Hodge
  theory of the Hecke category}

\author{Geordie Williamson}
\begin{document}

\begin{abstract}
Ideas from Hodge theory have found important applications in
representation theory. We give a survey of joint work with Ben Elias
which uncovers Hodge theoretic structure in the
Hecke category (``Soergel bimodules''). We also outline similarities
and differences to other combinatorial Hodge theories.
\end{abstract}



\maketitle

\section{Introduction}

In representation theory the setting is often highly structured. Thus it is surprising that several central theorems and
conjectures may be interpreted as stating that a given situation
behaves as if it were generic.\footnote{I learnt this point of view
  from W. Soergel and P. Fiebig.} (An example: a generic bilinear form on
a vector space is non-degenerate, however establishing that a specific
form is non-degenerate might be very difficult.)

The notion of a Hodge structure first arose in complex algebraic
geometry. 
The existence of Hodge structures is the extra data that
distinguishes the cohomology groups of algebraic varieties (or
K\"ahler manifolds) from those of
general spaces. Over the last twenty years it has been discovered that
some kind of combinatorial Hodge structure exists in at least four
other situations (polytopes, Coxeter groups, matroids, tropical geometry). The
existence of these structures has lead to the solutions of several conjectures
which have elementary formulations, but no elementary proofs. (An
example of such an elementary question: how many faces of each
dimension can a polytope have?) It is an
interesting question as to whether there is a unifying framework for
these various ``Hodge theories''.

In this paper we survey these Hodge theories and then concentrate on
Soergel bimodules, which gives rise to Hodge structures in three
distinct ways
(global, relative and local). Hodge structures are powerful enough
to deduce the genericity statements alluded
to above, and thus to prove several positivity conjectures
about the Kazhdan-Lusztig basis. They also provide new
algebraic proofs of difficult theorems in representation theory  (Kazhdan-Lusztig
conjecture, Jantzen conjecture).

The structure of this survey is as follows. In \S \ref{cht} we define
what we mean by a combinatorial Hodge theory and in \S \ref{ex} we
give some examples. In \S \ref{hecke} we recall the Hecke algebra, its
Kazhdan-Lusztig basis and its categorification by Soergel
bimodules. Finally, in \S \ref{htsb} we survey the
Hodge structures arising from Soergel bimodules.

%\tableofcontents

\section{Combinatorial Hodge theory} \label{cht}

Here we outline what we mean by combinatorial Hodge
theory. So far most of the structures that show up outside complex
algebraic geometry are of a very simple form (they are of
``Hodge-Tate type''). This means that much of the linear algebra
simplifies greatly.

\subsection{Lefschetz data} \label{lef}
Let $\Lambda \subset \RM$ denote a subfield. Examples to keep
in mind are $\Lambda = \QM, \RM$ or some finite extension of $\QM$. By
\emph{Lefschetz data} we mean the following:
\begin{enumerate}
\item A finite dimensional graded $\Lambda$-vector space $H = \bigoplus_{i \in
    \ZM} H^i$ which vanishes in either even or odd degree.
\item A non-degenerate graded symmetric bilinear form
\[
\langle -, - \rangle : H \times H \to \Lambda
\]
(``graded'' means that  $\langle H^i, H^j \rangle = 0$ if $i \ne -j$).
\item A graded vector space $V$ concentrated in degree 2 together with
  a map of graded algebras
\[
a : S^\bullet(V) \to \End(H)
\]  
where $S^\bullet(V)$ denotes the graded symmetric algebra on $V$
(i.e. we are given commuting endomorphisms $a(\gamma) : H^{\bullet}
\to H^{\bullet + 2}$ depending linearly on $\gamma \in V$).
Via $a$ we may view $H$ as a graded
$S^\bullet(V)$-module. From now on we forget $a$ and, for any $p \in
S^\bullet(V)$, simply write
\[
p \cdot h := a(p)(h).
\]
We require that the action is compatible with $\langle -, -
\rangle$ in the sense that
\[
\langle p \cdot h, h' \rangle = \langle h, p \cdot h' \rangle \quad \text{for all $p \in S^\bullet(V)$ and $h, h' \in H$.}
\]
\item An open convex cone $V_\ample \subset V$ (``cone'' means
  that $V_\ample$ is closed under multiplication by $\Lambda_{>0} :=
  \RM_{> 0} \cap \Lambda.$) 
\end{enumerate}

%\begin{definition}
We say that $\gamma \in V$ \emph{satisfies hard Lefschetz} if for all
$i \ge 0$ the action by $\gamma^i$ yields an isomorphism
\[
\gamma^i : H^{-i} \simto H^i.
\]
We say that Lefschetz data as above \emph{satisfies hard Lefschetz} if all $\gamma
\in V_\ample$ satisfy hard Lefschetz.
%\end{definition}

Let us fix $\gamma \in V$ which satisfies hard Lefschetz. We define the
\emph{$\gamma$-primitive subspace} (or simply \emph{primitive
  subspace}) to be
\[
P_\gamma^{-i} := \ker( \gamma^{i+1} : H^{-i} \to H^{i+2}) \subset H^{-i}.
\]
For any such $\gamma$ the map induced by the inclusions $P_\gamma^{-i}
\hookrightarrow H^{-i}$ gives an isomorphism of $\Lambda[\gamma]$-modules (the \emph{primitive
  decomposition}):
\[
\bigoplus_{i \ge 0} \Lambda[\gamma]/(\gamma^{i+1}) \otimes_\Lambda
P^{-i}_\gamma \simto H.
\]

\begin{remark} \label{rem:sl2}
Let us fix Lefschetz data $H, V$ etc.~as above.  Consider the Lie algebra $\mathfrak{sl}_2 := \Lambda f \oplus \Lambda h
  \oplus \Lambda e$ over $\Lambda$ with
\[
f = \left ( \begin{matrix} 0 & 0 \\ 1 &
      0 \end{matrix} \right), \quad
h = \left ( \begin{matrix} 1 & 0 \\ 0 & -1
       \end{matrix} \right) \quad \text{and} \quad 
e = \left ( \begin{matrix} 0 & 1 \\ 0 &
      0 \end{matrix} \right).
\]
Then $\gamma \in V$ satisfies hard Lefschetz if and only if there is
a (necessarily unique) $\mathfrak{sl}_2$-action on $H$ with $e(x) = \gamma \cdot x$ and $h(x) = j x$ for
all $x \in H^j$ (see e.g.~\cite[\S 1]{Patimo}). With respect to this action, the primitive decomposition is the isotypic
decomposition and the primitive subspaces
are the lowest weight spaces (i.e. $P_\gamma^{-i} = \ker (f : H^{-i}
\to H^{-i-2})$).
\end{remark}

For any $\gamma \in V$ and $i \ge 0$ we can consider the symmetric form
\begin{align*}
(-, -)_{\gamma}^{-i} : H^{-i} \times H^{-i} &\to \Lambda \\
(h,h') &\mapsto \langle h, \gamma^i \cdot h' \rangle  
\end{align*}
It is non-degenerate if and only if $\gamma$ satisfies hard Lefschetz on $H$.
Let $m$ denote the minimal non-zero degree in $H$. We say that $\gamma \in
V$ \emph{satisfies the Hodge-Riemann bilinear
  relations} if for all $i \ge 0$ the restriction of
$(-,-)_\gamma^{-i}$ to $P_{\gamma}^{-i} \subset H^{-i}$ is
$(-1)^{(-i - m)/2}$-definite.\footnote{The  exponent $(-i -
    m)/2$ is always an integer if $H^{-i} \ne 0$ by our assumption
    that $H$ vanishes in even or odd degrees.} In other words, we
  require that
  \begin{itemize}
  \item $(-,-)_\gamma^{m}$ is positive definite on $P_\gamma^m = H^m$,
  \item $(-,-)_\gamma^{m+2}$ is negative definite on $P_\gamma^{m+2} \subset H^{m+2}$,
  \item $(-,-)_\gamma^{m+4}$ is positive definite on $P_\gamma^{m+4} \subset
    H^{m+4}$, etc.
  \end{itemize}
If $\gamma$ satisfies the
Hodge-Riemann bilinear relations then it satisfies hard Lefschetz. We
say that Lefschetz data as above \emph{satisfies the Hodge-Riemann
  bilinear relations} if all $\gamma \in V_\ample$ do.

\begin{remark}
The set of $\gamma \in V$ which satisfy hard Lefschetz is Zariski
open and stable under multiplication by $\Lambda^\times$. The set of
$\gamma \in V$ which satisfy the Hodge-Riemann bilinear relations is
an open semi-algebraic set stable under multiplication by
$\Lambda_{>0}$.
\end{remark}

\subsection{The case of a graded algebra} \label{alg}
Many examples of Lefschetz data arise from Frobenius algebras. Consider
\begin{enumerate}
\item A positively graded commutative $\Lambda$-algebra $A =
  \bigoplus_{i \ge 0} A^i$ which
  vanishes in odd degree.
\item A degree zero map (``Frobenius form'')
\[
\tr : A \to \Lambda(-2m)
\] 
which does not vanish on any non-zero ideal of $A$. (Here $\Lambda(-2m)$
denotes the one-dimensional graded vector space concentrated in degree
$2m$, for some $m \ge 0$.)
\item An open convex cone $A^2_\ample \subset A^2$.
\end{enumerate}

To this we may associate a Lefschetz datum as follows:
\begin{enumerate}
\item We set $H := A(m)$, i.e. $H = \bigoplus H^i$ where $H^i := A^{i+m}$.
\item We define $\langle a, a' \rangle := \tr(aa')$ 
  (non-degenerate by assumption (2)).
\item $V_\ample := A^2_\ample \subset V := A^2$ with action induced by
  the multiplication in $A$.
\end{enumerate}

\section{Some examples} In this section we give some examples of
Lefschetz data. \label{ex}

\subsection{Classical Hodge theory} Let $X$ be a connected smooth complex projective
variety (or more generally a compact K\"ahler manifold) of (complex) dimension $m$. Let $A := H^*(X, \QM)$ denote its
(singular) cohomology ring. It is equipped with the trace $\tr : A^{2m} \to
\QM$ given by pairing with the fundamental class of $X$ and the
corresponding Poincar\'e form $\langle \alpha, \beta \rangle = \tr( \alpha \cup \beta)$.
Let  $A^2_\ample \subset A^2 = H^2(X, \QM)$ denote the ample cone (the convex
hull of all $\QM_{>0}$-multiples of Chern classes of ample line bundles).

Assume that the Hodge decomposition of $H^*(X, \CM)$ involves only
type $(p,p)$.\footnote{We ignore here the more general form of the
Hodge-Riemann bilinear relations (see e.g. \cite{Voisin}) when Hodge types other than $(p,p)$ are
present.} By classical Hodge theory (see
e.g. \cite{Voisin}) the corresponding Lefschetz datum (see \S \ref{alg}) satisfies hard
Lefschetz and the Hodge-Riemann bilinear relations.

If $X$ is still projective but no longer smooth we can instead consider its
intersection cohomology $H := IH^*(X, \QM)$. This is a graded
$H^*(X)$-module concentrated in degrees between $-m$ and $m$, and is
equipped with a non-degenerate intersection form $\langle -,
-\rangle$. Suppose as before that the Hodge decomposition of $H \otimes_{\QM} \CM$ only
involves type $(p,p)$. Then the Lefschetz datum $(H, \langle -,
-\rangle, A_\ample)$ satisfies hard Lefschetz and the Hodge-Riemann
bilinear relations \cite{BBD, dCM2}.

In examples coming from algebraic geometry, all Lefschetz data
are naturally defined over $\QM$.

\subsection{Polytopes}  Recall that a polytope $P \subset \RM^n$ is the
convex hull of finitely many points in $\RM^n$. A polytope is
simplicial if all of its maximal dimensional faces are simplices. 
Any $n$-dimensional polytope determines a collection $\Delta$ of polyhedral cones (a
``fan'') by moving $P$ so that its interior contains the origin and considering
the cones spanned by all faces $F \ne P$ 
of $P$. To $P$ we may associate $A(P)$, the algebra of
functions on $\RM^n$ which are piecewise polynomial on all cones in $\Delta$, and its
quotient $\Pi(P)$ by the action of polynomials of positive degree on
$\RM^n$. The algebras $A(P)$ and $\Pi(P)$ are graded by degree. Inside
$\Pi(P)^1$ we can consider the convex cone of
$\Pi(P)^1_\ample$ given by (the image of) strictly convex piecewise
linear functions on $\Delta$. If $P$ is simplicial then there exists a
trace map $\tr : \Pi(P)^n \to
\RM$ \cite{Brion}. (In this case $\Pi(P)$ is isomorphic to the \emph{polytope
  algebra} of McMullen \cite{McM1}, and under this isomorphism $\tr$
is given by  volume.\footnote{More precisely, we should tensor the
degree zero part of the polytope algebra by the real numbers to obtain
such an isomorphism, see \cite{McM1}.})

If $P$ is simplicial then (after doubling degrees and applying the
recipe of
\S \ref{alg}) the associated Lefschetz data satisfies hard Lefschetz and the Hodge-Riemann bilinear
relations. The hard Lefschetz theorem can be used to deduce the
necessity of McMullen's conditions on the number of faces of $P$ of
each dimension (``face numbers'')
\cite{McM0}. The Hodge-Riemann relations imply generalisations of the
Aleksandrov-Fenchel inequalities in convex geometry. The hard
Lefschetz theorem and Hodge-Riemann bilinear relations were proved by
Stanley by ``wiggling the vertices'' of $P$ so that they become rational
numbers, and then identifying $\Pi(P)$ with the cohomology ring of a
rationally smooth projective toric variety \cite{Stanley}. A direct
proof (entirely within the world of convexity) was given by McMullen
\cite{McM2}.

If $P$ is no longer simplicial then Bressler and Lunts
\cite{BrLu} described a
combinatorial recipe to produce a graded $\Pi(P)$-module $H(P)$, which agrees
with the intersection cohomology of a projective (but no longer
necessarily rationally smooth) toric variety if the vertices of $P$ are rational
numbers (similar results were obtained by Barthel, Brasselet, Fieseler
and Kaup \cite{BBFK}). The graded module $H(P)$ is
equipped with a symmetric non-degenerate form $\langle -, -\rangle :
H(P) \times H(P) \to \RM$. It is a theorem of Karu \cite{Karu} (see
also \cite{Braden}) that the Lefschetz data $(H(P), \langle -,
-\rangle, \Pi(P)^1_\ample)$ satisfies hard Lefschetz and the
Hodge-Riemann bilinear relations. The hard Lefschetz theorem is
used to deduce the unimodality of Stanley's generalised $h$-vector,
however here it is unclear if these inequalities have any significance
for face numbers.

The (initially counterintuitive) fact that there exist polytopes whose
combinatorial type is not that of any polytope with rational
vertices (see e.g. \cite{Ziegler}) means that the trick of ``wiggling
vertices'' cannot be used to
deduce Karu's theorem from the hard
Lefschetz theorem and Hodge-Riemann relations for algebraic varieties.

Note that, if $P$ has vertices defined over some subfield
$\Lambda \subset \RM$ then $\Pi(P)$ and $H(P)$ are defined over
$\Lambda$. Hence non-rational polytopes give many examples of
Lefschetz data with no natural $\QM$-structure.

\subsection{Coxeter groups} 
This is the subject of this survey. The theory of Soergel bimodules
yields Lefschetz data satisfying hard Lefschetz and the Hodge-Riemann
bilinear relations in three distinct ways (which we refer to as 
\emph{global}, \emph{relative} and \emph{local}). The
simplest interesting example of a vector space underlying Lefschetz
data that arises from this theory is the
\emph{coinvariant algebra} (see Remark \ref{w0})
\[
C := R / \langle R^W_+ \rangle.
\]
Here $W \subset GL(V)$ denotes a finite reflection group;
$R$ denotes the polynomial functions on $V$, graded so that
  $\deg V^* = 2$;
$\langle R^W_+ \rangle$ denotes the ideal generated by $W$-invariants
which are homogenous of strictly
positive degree.

Examples arising from Soergel bimodules are defined over finite
extensions of $\QM$ obtained by adjoining algebraic integers of the
form $2\cos(\pi/m_{st})$ where $m_{st} \in \ZM_{\ge 2}$. Thus there
are many examples with no natural $\QM$-structure (and probably no
$\QM$-structure at all). Thus it seems unlikely that they arise from
complex algebraic varieties in any straightforward way.

\subsection{Matroids}
Let $M$ be a matroid of rank $r + 1$ on a set of cardinality
$n+1$. (For example, if $k$ is a field and if $\{ v_0, \dots,
v_n\} \subset V$ is a spanning set of vectors in an $r+1$-dimensional
$k$-vector space then we obtain a matroid which tells us which subsets of $\{ v_0, \dots,
v_n\}$ are linearly independent, of maximal size and fixed
rank, etc. Matroids arising from such arrangements are said to be
\emph{realisable} over $k$.) To $M$ one may associate its \emph{Chow
  ring} $A^*(M)$, a graded commutative $\RM$-algebra. Under mild assumptions
(``loopless'') on
$M$, the Chow ring
is equipped with an isomorphism $\deg : A^r(M) \simto \QM$ and a cone
$A^1(M)_{\textrm{conv}} \subset A^1(M)$ of ``strictly
convex functions'' (see \cite{AHK} and the references therein).

In \cite{AHK}, it is proved (after doubling degrees and applying
the recipe of \S \ref{alg}) that the associated
Lefschetz datum satisfies hard Lefschetz and the Hodge-Riemann
bilinear relation. From the Hodge-Riemann relations the authors deduce the log concavity of
the absolute value of the coefficients of the characteristic
polynomial of $M$, an old conjecture in matroid theory. This
generalises a similar log concavity property for the chromatic
polynomial of a graph \cite{Huh}.

The Chow ring $A^*(M)$ is defined over $\QM$. In \cite[\S 5.4]{AHK} it is proved that $A^r(M)$ is naturally isomorphic
to the Chow ring of a rationally smooth projective variety over $k$ if
and only if the matroid is realisable over $k$. (However recall that ``most''
matroids are probably not realisable over \emph{any} field.) If $k = \CM$ then the
hard Lefschetz and Hodge-Riemann relations for $A^*(M)$ may be deduced from
classical Hodge theory; for general $k$ they are
related to Grothendieck's standard conjectures.

\begin{remark}
The reader is referred to \cite{Tropical} for another remarkable
example (of quite a different flavour to those discussed above) of
combinatorial Hodge theory, this time arising from tropical geometry.
\end{remark}

\section{The Hecke category} \label{hecke}

\subsection{Coxeter systems}

Let $(W,S)$ denote a Coxeter system. That is, $W$ is a group together
with a distinguished finite generating set $S \subset W$ of
\emph{simple reflections} such that $W$ admits a \emph{Coxeter
presentation}
\begin{align*}
  W = \langle s \in S \;|\; (st)^{m_{st}} = \textrm{id} \text{ for all $s,
    t \in S$} \rangle
\end{align*}
for certain $m_{st} \in \ZM_{\ge 0}$ such that $m_{ss} = 1$ for all $s
\in S$ and $m_{st} \in \{2, 3, 4, \dots \} \cup \{ \infty \}$ for $s
\ne t$ ($m_{st} = \infty$ means that we do not impose
any relation on $st$). By definition the \emph{reflections} is the subset $T \subset
W$ consisting of the conjugates of $S$:
\[
T := \bigcup_{w \in W} wSw^{-1}.
\]

\begin{remark}
Coxeter groups are so called because of Coxeter's theorem that any
finite reflection group (i.e. finite subgroup of linear transformations of
a real vector space generated by reflections) may be given a Coxeter
presentation.
\end{remark}

For any $x \in W$ a reduced expression is an expression $\un{x} =
(s_1, \dots, s_m)$ for $x$ in the generating set $S$ (i.e. $x = s_1
\dots s_m$ with all $s_i \in S$) which is of minimal length. We
denote by $\ell : W \to \ZM_{\ge 0}$ the length function on $W$
(i.e. $\ell(x)$ is the length of a reduced expression for $x$). Let $\le$ denote the Bruhat order on
$W$.

Central to the theory of Soergel bimodules is a certain reflection\footnote{i.e. reflections act as
  reflections} 
representation of $W$. To this
end consider a finite dimensional real vector space $\hg$ together
with subsets $\{ \alpha^\vee_s \}_{s \in S} \subset \hg$ and $\{
\alpha_s \}_{s \in S} \subset \hg^*$ of \emph{coroots} and
\emph{roots} satisfying the following two conditions:
\begin{enumerate}
\item the subsets $\{ \alpha^\vee_s \}_{s \in S} \subset \hg$ and $\{
\alpha_s \}_{s \in S} \subset \hg^*$ are linearly independent;
\item under the natural pairing $\hg^* \times \hg \to \RM$ we have
  \begin{equation}
    \label{eq:cos}
\langle \alpha_s, \alpha_t^\vee \rangle = -2 \cos(\pi/m_{st}).    
  \end{equation}
\end{enumerate}
Then the assignment $s
\mapsto \phi^\vee_s \in GL(\hg)$ (resp. $s \mapsto \phi_s \in GL(\hg^*)$) where
  \begin{align*}
    \phi^\vee_s(v) := v - \langle \alpha_s , v\rangle \alpha_s^\vee
\qquad \quad ( \text{resp.} \quad
    \phi_s(\lambda) := \lambda - \langle \lambda , \alpha_s^\vee \rangle \alpha_s)
  \end{align*}
defines a representation of $W$ on $\hg$ (resp. $\hg^*$).

\begin{remark} Some remarks are in order:
\begin{enumerate}
\item The above definition mimics the action of the Weyl group on a
  Cartan subalgebra of a
  complex semi-simple Lie algebra. The only difference to the Weyl
  group case is  that we have scaled all roots so as to have the same
  length.
\item Central to the study of general Coxeter systems is a ``geometric''
  representation first defined by Tits. One sets $\hg := \bigoplus_{s
    \in S} \RM \alpha_s^\vee$ and defines $\alpha_s \in
  \hg^*$ by equation \eqref{eq:cos} above. With this definition $\{ \alpha_s
  \}_s \in \hg^*$ will be linearly independent if and only if
$W$ is finite. Thus the representation considered above usually does
not agree with Tits' representation. The need to ``enlarge'' so that both
roots and coroots are linearly independent is important in
the theory of Kac-Moody Lie algebras \cite{Kac}. It is also aesthetically
pleasing to have $\hg$ and $\hg^*$ play symmetric roles.
\end{enumerate}  
\end{remark}

By our assumptions above the intersections of half-spaces
\begin{align*}
  \hg_\reg^+ := \bigcap_{s \in S} \{ v \in \hg \; | \; \langle \alpha_s, v \rangle > 0 \}
\quad \text{and} \quad
  \hg^{*+}_\reg := \bigcap_{s \in S} \{ \lambda \in \hg \; | \; \langle \lambda, \alpha_s^\vee \rangle > 0 \}
\end{align*}
are non-empty. Borrowing terminology from Lie theory we refer to
elements in either set as \emph{dominant regular}.

\subsection{The Hecke algebra and Kazhdan-Lusztig basis}
The Hecke algebra is a deformation of the group algebra of a Coxeter
group which plays a fundamental role in Lie theory. It is a free
$\ZM[v^{\pm 1}]$-algebra $H$ with basis $\{ h_x \}_{x \in W}$ and multiplication
determined by the rules (for $s \in S$ and $x \in W$)
\[
h_sh_x = \begin{cases} h_{sx} & \text{if $sx > x$,} \\
(v^{-1} - v)h_x + h_{sx} & \text{if $sx < x$}. \end{cases}
\]
The basis $\{ h_x \; | \; x \in W \}$ is the \emph{standard basis} of
$H$. The algebra $H$ posesses an involution $h \mapsto \overline{h}$
determined by $v \mapsto v^{-1}$ and $h_x \mapsto h_{x^{-1}}^{-1}$.

The algebra $H$ possesses a remarkable basis, discovered by
Kazhdan and Lusztig \cite{KL1}. The Kazhdan-Lusztig basis is the unique basis $\{ b_x \}$ for $H$ such that:
\begin{enumerate}
\item $\overline{b_x} = b_x$ (``self-duality'');
\item $b_x \in h_x + \sum_{y < x} v\ZM[v] h_y$.
\end{enumerate}
For example $b_s = h_s + vh_{\id}$. If we write
\[
b_x = \sum_{y \in W} p_{y,x} h_y
\]
the polynomials $p_{y, x}$ are \emph{Kazhdan-Lusztig polynomials}.

The Kazhdan-Lusztig basis appears to satisfy numerous (a priori very
mysterious) positivity
properties. Here is an incomplete list:
\begin{enumerate}
\item \emph{Positivity of Kazhdan-Lusztig polynomials:} 
\begin{equation}
  \label{pos1}
p_{y,x} \in \ZM_{\ge 0}[v].
\end{equation}
\item \emph{Positivity of inverse Kazhdan-Lusztig polynomials:} 
If we write
\[
h_x = \sum (-1)^{\ell(x) - \ell(y)} g_{y,x} b_y
\]
then
\begin{equation}
  \label{pos2}
g_{y,x} \in \ZM_{\ge 0}[v].
\end{equation}
\item \emph{Positivity of structure constants:} 
If we write
\[
b_xb_y = \sum \mu_{x,y}^z b_z
\]
then
\begin{equation}
  \label{pos3}
\mu_{x,y}^z
  \in \ZM_{\ge 0}[v^{\pm 1}].
\end{equation}
\item \emph{Unimodality of structure constants:} 
If we set
\begin{equation}
  \label{eq:gauss}
[m] := \frac{v^m-v^{-m}}{v-v^{-1}} = v^{-m+1} + v^{-m+3} + \dots +
v^{m-3} + v^{m-1}  
\end{equation}
and, for all $x, y, z \in W$, write
\[
\mu_{x,y}^z = \sum_{m \ge 1} a_{x,y}^{z,m}[m]
\]
(this is possible since $\overline{\mu_{x,y}^z} = \mu_{x,y}^z$), then
\begin{equation}
  \label{pos4}
a_{x,y}^{z,m} \in \ZM_{\ge 0}.
\end{equation}
(In other words, $\mu_{x,y}^z$ is the character of a finite
dimensional $\mathfrak{sl}_2(\CM)$-module.)
\end{enumerate}

\begin{remark}
As the reader has probably already noticed, \eqref{pos4} is strictly
stronger than \eqref{pos3}. We separate them here because (as we will
see) the explanation for \eqref{pos4} lies deeper than that of \eqref{pos3}.
\end{remark}

One of the main purposes of this note is to state theorems of a Hodge
theoretic nature about Soergel bimodules which imply:

\begin{theorem}[\cite{EWHodge} and \cite{EWrel}] Positivity properties
  \eqref{pos1}, \eqref{pos2}, \eqref{pos3} and \eqref{pos4} hold for
  any Coxeter system.
\end{theorem}

\begin{remark} That \eqref{pos1} holds for Weyl and affine Weyl groups
  is due to Kazhdan-Lusztig \cite{KL2}. That \eqref{pos2} is true for
  affine and Weyl groups is proved in \cite{Springer} (the proof
uses ideas of Springer, MacPherson and
  Brylinski). Since the construction of Kac-Moody flag varieties for
  any generalised Cartan matrix it was understood that the arguments
  of \cite{KL2,Springer} can be used to show that  \eqref{pos1}
  and \eqref{pos2} hold for any ``crystallographic'' Coxeter group
  (i.e. $m_{st} \in \{ 2, 3,4,6,\infty \}$ for all $s \ne t$).
\end{remark}

\begin{remark}
  There are almost no cases so far where there exists combinatorial
  (i.e. only using combinatorics of Coxeter groups) 
  proofs of the above positivity properties. Notable exceptions
  include the case of ``Grassmannian'' permutations in the symmetric
  group \cite{LS}, and the case of ``universal'' Coxeter
  systems \cite{Dyer}. Let us repeat Bernstein's opinion on the
  subject of combinatorial formulas
  \cite{Bernstein}: ``In some cases one can get explicit formulas for
  [Kazhdan-Lusztig polynomials]. For instance, one can calculate
  intersection cohomology for Schubert varieties on usual
  Grassmannians (see Lascoux and Sch\"utzenberger). But Zelevinsky
  showed that in this case it is possible to construct small
  resolutions of singularities. I would say that if you can compute a
  polynomial $P$ for intersection cohomologies in some case without a
  computer, then probably there is a small resolution which gives it.''
\end{remark}

\begin{remark}
The reader is referred to \cite{DL} for a list of further
  positivity properties and to \cite{Gobet} for some recent results.
\end{remark}

\subsection{Bimodules}

Let $R$ denote the regular functions on $\hg$. (After choosing a
basis $x_1, \dots, x_n$ for $\hg^*$, $R$ is simply the polynomial ring
$\RM[x_1, \dots, x_n]$.) Let $\bim$ denote the category of graded
$R$-bimodules, with morphisms those $R$-bimodule maps which preserve
the grading (i.e. are of ``degree zero''). Given bimodules $M, N \in
\bim$ we let $\Hom(M,N)$ denote the vector space of bimodule
homomorphisms. We denote the tensor product of $M$ and $N$ by juxtaposition:
\[
MN := M \otimes_R N \in \bim.
\]
This gives $\bim$ the structure of monoidal category with unit
$R$. For $M = \bigoplus_{i \in \ZM} M^i \in \bim$ and $m \in \ZM$ we
denote the shifted bimodule by $M(m)$
\[
M(m)^i := M^{m+i}.
\]
(So $M(m)$ is the same bimodule with grading shifted ``down'' by $m$.)
Given bimodules $M, N \in \bim$ we denote by $\Hom^\bullet(M,N)$ the
graded $R$-bimodule consisting of morphisms of all degrees:
\[
\Hom^\bullet(M,N) := \bigoplus_{m \in \ZM}\Hom(M, N(m)).
\]
(We emphasise that this is simply notation: it does not refer to the
homomorphisms in any category we consider in this paper.)

Let $\bimf$ denote the full subcategory of $\bim$ consisting of
bimodules which are finitely generated both as left and right
$R$-modules. Then $\bimf$ is a monoidal subcategory of $\bim$. If $M, N \in \bimf$ then $\Hom(M,N)$ is
finite dimensional. It follows that the additive category $\bimf$ is Krull-Schmidt: any
object admits a decomposition into indecomposable objects and an object is
indecomposable if and only if its endomorphism ring is local; it
follows that the Krull-Schmidt theorem holds in $\bimf$.

\subsection{Soergel bimodules} Recall that $R$ is the polynomial ring
of regular functions on $\hg$. Thus our Coxeter group $W$ acts on $R$ (by
functoriality). Given $x \in W$ we denote by $R^x$ the subring of
invariants under $x$. Similarly, $R^{x_1, x_2}$ denotes the invariants
under the subgroup generated by $x_1, x_2$ etc. For $s \in S$ we consider the bimodule
\[
B_s := R \otimes_{R^s} R (1).
\]
It is easy to see that $B_s$ belongs to $\bimf$. We denote by $\HC$
the smallest strict, graded, additive, monoidal and
Karoubian\footnote{i.e. $\HC$ is closed under taking direct
  summands: if $M \oplus M' = N \in \HC$ then $M, M' \in
  \HC$} subcategory of
$\bimf$ containing $B_s$ for all $s \in S$.
In formulas:
\[
\HC := \langle B_s \;|\; s \in S \rangle_{\cong, (\pm 1), \oplus, \otimes,
  \textrm{Kar}} \subset \bimf.  
\]
By the Krull-Schmidt theorem, the indecomposable
objects of $\HC$ are the indecomposable direct summands of the \emph{Bott-Samelson
  bimodules}
\[
BS(\un{x}) = B_sB_t \dots B_u(m) \in \HC
\]
for all expressions $\un{x} = (s, t, \dots, u)$ in the simple
reflections and all $m \in \ZM$. Bimodules belonging to $\HC$ are called
\emph{Soergel bimodules}. We call the category $\HC$
the category of Soergel bimodules or the \emph{Hecke category} (see Remark
\ref{rem:hc}).

\begin{remark}
  The reader is warned that $\HC$ is additive but never abelian.
\end{remark}

It is a remarkable theorem of Soergel that $\HC$ categorifies the
Hecke algebra. Let $[\HC]$ denote the split Grothendieck group of
$\HC$: it is the free abelian group generated by symbols $[M]$ for $M
\in \HC$ modulo the relation $[M] = [M'] + [M'']$ if $M \cong M'
\oplus M''$. We view $[\HC]$ as an $\ZM[v^{\pm 1}]$-algebra via:
\begin{align*}
  [M][N] &:= [M \otimes N'], \\
  v[M] &:= [M(1)].
\end{align*}

\begin{theorem}[``Soergel's categorification theorem'', \cite{SoeB}]
  There exists an (obviously unique) isomorphism of $\ZM[v^{\pm
    1}]$-algebras:
  \begin{align}
  \phi:  H  \simto [\HC] : b_s \mapsto [B_s].
\end{align}
\end{theorem}

To establish the existence of $\phi$ it is enough to verify certain
isomorphisms among tensor products of $B_s$ and $B_t$ for pairs
$s, t \in S$ (see the the examples
below). In the words of Soergel \cite{SoeICM}: ``This is a bit
tricky, but not deep''.  To show that $\phi$ is an isomorphism we must control the size of
$[\HC]$, which comes as a corollary of the following theorem:

\begin{theorem}[\cite{SoeB}]
For each $x \in W$ there exists a unique indecomposable bimodule $B_x$ (well-defined
up to isomorphism) which occurs as a summand of the Bott-Samelson bimodule $BS(\un{x})$ for any
reduced expression $\un{x}$ of $x$, and does not occur as a summand of
$BS(\un{y})$ for any shorter expression $\un{y}$. The set
\[
\{ B_x(m) \; | \; x \in W, m \in \ZM \}
\]
gives a set of representatives for the indecomposable Soergel
bimodules up to isomorphism.
\end{theorem}

It follows that the classes of the bimodules $B_x$ give a basis for $[\HC]$:
\begin{equation}
  \label{eq:basis}
[\HC]  = \bigoplus_{x \in W} \ZM[v^{\pm 1}] [B_x].  
\end{equation}
Once one has established the existence of a homomorphism $\phi : H \to
[\HC]$, 
Soergel's categorification theorem is easily deduced from \eqref{eq:basis}.

\begin{example}
  \begin{enumerate}
  \item Suppose that $S = \{ s \}$ so that $W = \ZM/2\ZM$ acting on $R =
    \RM[\alpha]$ via $ s(\alpha) = -\alpha$, and $R^s =
    \RM[\alpha^2]$. In this case one may easily calculate that $b_s =
    h_s + v$ satisfies
\[
b_sb_s = (v+v^{-1})b_s.
\]
In this
    case Soergel's categorification theorem follows easily from the
    isomorphism
    \begin{align}
\nonumber
B_sB_s & = R \otimes_{R^s} R \otimes_{R^s} R(2) \\ \nonumber
& \cong R \otimes_{R^s} ( R^s \oplus R^s(-2)) \otimes_{R^s} R(2) \\ 
& = B_s(1) \oplus B_s(-1).            \label{eq:Bs}
    \end{align}
where for the middle step we use that $R = R^s \oplus \alpha R^s$ as
an $R^s$-bimodule: ``any
polynomial can be written as the sum of an even and an odd
polynomial''. Thus the indecomposable Soergel bimodules are (up to
shift) given by
\[
\{ R := B_\id, B_s \}.
\]
\item Suppose that $W = \langle s, t \rangle$ with $m_{st} = 3$ so
  that $W$ is isomorphic to the symmetric group on three letters. In
  this case one may check easily that
  \begin{align*}
    b_sb_t = b_{st},  \quad b_sb_tb_s = b_{sts} +
    b_s\quad \text{and} \quad b_{sts}b_s = (v+v^{-1})b_{sts}.
  \end{align*}
(and similarly with $s$ and $t$ interchanged). In this case Soergel's
categorification amounts to the isomorphism \eqref{eq:Bs} for $B_s$ and
$B_t$ as well as the statements (which can be checked by hand):
\begin{gather*}
  B_{st} := B_sB_t \quad \text{is indecomposable as an $R$-bimodule,} \\
B_sB_tB_s \cong B_{sts} \oplus B_s \quad \text{with $B_{sts} := R
  \otimes_{R^{s,t}} R(3)$,} \\
B_{sts}B_s \cong B_{sts}(1) \oplus B_{sts}(-1).
\end{gather*}
(and similarly with $s$ and $t$ interchanged). Indeed
once one knows \eqref{eq:Bs} and the above three statements then it is easy
to show by induction that any indecomposable summand of any tensor
product of $B_s$ and $B_t$ is isomorphic to a shift of one of the
indecomposable bimodules
\[
\{ B_\id := R, B_s, B_t, B_{st}, B_{ts}, B_{sts} \}.
\]
  \end{enumerate}
\end{example}

\begin{remark}
  The above examples are deceptive: in general there is no explicit
  description of the indecomposable Soergel bimodules (just as there
  is no explicit description of the Kazhdan-Lusztig
  polynomials). Soergel's proof \cite{SoeB} is highly
  non-constructive. It formally ressembles the proof giving uniqueness of tilting
  objects in highest weight categories \cite{Ringel, Donkin}.
\end{remark}

In \cite{SoeB}, Soergel also constructs an explicit inverse to $\phi$ (``the character of a
Soergel bimodule''):
\[
\ch : [\HC] \stackrel{\sim}{\to} H.
\]
Its definition (which we do not give here) involves taking the graded
rank of the subquotients of certain filtrations. It is manifestly
positive on bimodules, i.e.
\[
\ch(M) := \ch([M]) \in \bigoplus_{x \in W} \ZM_{\ge 0}[v^{\pm 1}] h_x
\quad \text{for  $M \in \HC$.}
\]
Thus $\{ \ch(B_x) \; | \; x \in W \}$
gives a basis for $H$ which is positive (i.e. belongs to $\bigoplus
\ZM_{\ge 0}[v^{\pm 1}] h_x$) and has structure constants in $\ZM_{\ge
  0}[v^{\pm 1}]$ (the structure constants give the
graded multiplicity of $B_z$ as a summand of $B_xB_y$).

\begin{conjecture}[Soergel \cite{SoeB}] For all $x \in W$, $\ch(B_x) =
  b_x$.
\label{conj:S}
\end{conjecture}

This conjecture is proved in \cite{EWHodge} using Hodge theoretic
ideas, as we will discuss below. It immediately implies the positivity
properties \eqref{pos1} and \eqref{pos3}. The method of proof also
establishes \eqref{pos2} (see \cite[Remark 6.10]{EWHodge}). A variant of
these methods establishes \eqref{pos4}.

\begin{remark} \label{rem:hc}
This conjecture was proved for Weyl groups and dihedral groups by Soergel
in
\cite{SoeHC}. It is proved for the Weyl group of any Kac-Moody Lie
algebra by H\"arterich in \cite{Ha} and for universal Coxeter systems
by Fiebig \cite{Fiebig}. Outside of dihedral groups and universal
Coxeter systems, these
proofs rely on the decomposition theorem: the bimodule $B_x$ is
realised as the equivariant intersection cohomology of a Schubert
variety.  Soergel's original
motivation for studying these bimodules came from attempts to
better understand the Kazhdan-Lusztig conjecture
\cite{KL1}; indeed his conjecture implies the Kazhdan-Lusztig
conjecture \cite{SoeJAMS}.
\end{remark}

\begin{remark} \label{rem:hc}
In geometric settings the category $\HC$ has other incarnations (as
semi-simple complexes on the flag variety, or as a variant of Harish-Chandra bimodules,
\dots). In \cite{EWSC}, a monoidal category is defined by
explicit diagrammatic generators and relations and it is proved that this category
is equivalent to $\HC$. The upshot is that the category of Soergel
bimodules is one incarnation of a more fundamental object, 
often referred to as the Hecke category. This point of view is
particularly useful when one wishes to study variants of the theory
over fields of positive characteristic, where the theory of Soergel
bimodules becomes unwieldy. The reader is referred to the
introduction to \cite{EWSC} for more on this point of
view. The Braden-MacPherson and Fiebig theory of sheaves on moment
graphs can be seen as another incarnation of the Hecke category
\cite{BM, Fiebig}.
\end{remark}

\section{Hodge theory of Soergel bimodules} \label{htsb}

\subsection{Global theory}

Recall the indecomposable Soergel bimodules $B_x$ introduced in the
previous section. In the global theory a central role is played by the
corresponding \emph{Soergel modules}, which are obtained by
quotienting by the action of positive degree polynomials on the right:
\begin{equation*}
  \overline{B}_x := B_x \otimes_R \RM.
\end{equation*}
Each $\overline{B}_x$ is a graded $R$-module, which is finite
dimensional over $\RM$ and vanishes in degree of parity different to
that of $\ell(x)$.

\begin{remark} \label{rem:smod}
  As already remarked above, in geometric settings
  $B_x$ may be obtained as the equivariant intersection cohomology of
  a Schubert variety. In such cases $\overline{B}_x$ is isomorphic to the ordinary
  (i.e. non-equivariant) intersection cohomology. See \cite{EWSurvey}
for a detailed discussion of Soergel modules and their connection to
intersection cohomology.
\end{remark}

Any indecomposable Soergel bimodule carries an \emph{intersection
  form}. This is a graded symmetric bilinear  form
\[
\langle -, -\rangle : B_x \times B_x \to R
\]
which is characterised (up to multiplication by $\RM_{>0}$) by the
following properties (see \cite[Lemmas 3.7 and 3.10]{EWHodge}):
\begin{enumerate}
\item For all $b, b' \in B_x$ and $r \in R$ we have:
  \begin{gather}
%    \langle b, b' \rangle = \langle b', b \rangle, \\
    \langle rb, b' \rangle = \langle b, rb' \rangle, \\
       \langle br, b' \rangle = \langle b, b' \rangle r  =\langle b,
       b'r \rangle. \label{eq:right}
\end{gather}
\item By \eqref{eq:right}, $\langle -, -\rangle$ descends to a form on the
corresponding Soergel module
\begin{equation*}
  \langle -, -\rangle_{\overline{B}_x}: \overline{B}_x \times
  \overline{B}_x \to R/R^{>0} = \RM
\end{equation*}
such that $\langle rm, m'\rangle = \langle m, rm'\rangle$ for all $r
\in R$ and $m, m' \in \overline{B}_x$. This form is non-degenerate.
\item If $b_\mini \in B_x$ denotes a non-zero element of degree
  $-\ell(x)$ (the minimal non-zero degree) then
\begin{equation}
  \label{eq:1}
\langle \lambda^{\ell(x)} b_\mini, b_\mini \rangle_{\overline{B}_x} > 0  
\end{equation}
for any (equivalently every) dominant regular $\lambda \in
\hg^*$.
\end{enumerate}

\begin{remark} The reader may well wonder where the intersection form
  on $B_x$ comes from. Recall that $B_x$ occurs as a direct
  summand of $BS(\un{x})$ for any reduced expression $\un{x}$. It is
  not difficult to equip $BS(\un{x})$ with a symmetric non-degenerate
  $R$-valued form (see \cite[\S 3.4]{EWHodge}). Somewhat miraculously,
  this form restricts to yield the desired form on $B_x$ (for any
  choice of embedding).
\end{remark}

\begin{remark}
  In geometric settings the intersection form may be identified with
  the (topological) intersection form on equivariant intersection
  cohomology.
\end{remark}

Thus $(\overline{B}_x, \langle -, -\rangle_{\overline{B}_x}, \hspos
\subset \hg^*)$ gives Lefschetz data (see \S \ref{lef}). The
main theorem of \cite{EWHodge} is the following:

\begin{theorem} \label{global}
  For any $x \in W$ and $\lambda \in \hspos$, the action of $\lambda$ on
  $\overline{B}_x$ satisfies hard Lefschetz and the Hodge-Riemann
  bilinear relations.
\end{theorem}

\begin{remark}
 This theorem is not new in
  geometric settings, where it may be deduced from the hard Lefschetz
  theorem and Hodge-Riemann bilinear relations for intersection
  cohomology.
\end{remark}

\begin{remark}
  We will not go into how Theorem \ref{global} is related to Soergel's
  conjecture (Conjecture \ref{conj:S}). In rough outline the
  strategy of proof in \cite{EWHodge} is to establish 
  Soergel's conjecture and Theorem \ref{global} ``at
  the same time''
  by an induction over the Bruhat order. Let us simply mention that ideas from de
  Cataldo and Migliorini's proof of the decomposition theorem play a
  key role \cite{dCM1, dCM2, WilBourbaki}. In particular, the definiteness provided by the
  Hodge-Riemann relations is central to the proof.
\end{remark}

\begin{remark} \label{w0}
  If $W$ is finite with longest element $w_0$ then $B_{w_0} = R
  \otimes_{R^W} R(\ell(w_0))$ ($R^W \subset R$ denotes the subalgebra of
  $W$-invariants) and hence
\[
\overline{B}_{w_0} = R/ \langle R^W_+ \rangle (\ell(w_0)).
\]
If $W$ is a finite Weyl group then it is well-known that this algebra
is (up to shift) isomorphic to the cohomology ring of the flag variety
and Theorem \ref{global} for $x = w_0$ follows from classical Hodge
theory. However if $W$ is not a Weyl group then there is no known
geometric proof of Theorem \ref{global} for $x = w_0$. (This case
seems like it should be much simpler than the general case.) The
surprisingly interesting example of
dihedral groups is discussed in detail in \cite[\S 6]{EWSurvey}.
\end{remark}

\begin{remark}
For a smooth projective variety $X$ Looijenga-Lunts \cite{LL} considered
the Lie subalgebra of endomorphisms of $H^*(X, \RM)$ generated all
copies of $\mathfrak{sl}_2(\RM)_\lambda$ associated to all ample classes
$\lambda \in H^2(X, \RM)$ (here $\mathfrak{sl}_2(\RM)_\lambda$ denotes
the copy of $\mathfrak{sl}_2(\RM)$ inside the endomorphisms of
$H^*(X,\RM)$ provided by Remark \ref{rem:sl2}). They show (using
the Hodge-Riemann bilinear relations) that one always obtains a
reductive Lie algebra in this way, and compute several examples. By
Theorem \ref{global} the definition of this Lie algebra also makes
sense for any Soergel module. Recently,
Patimo \cite{Patimo} has shown that for many $x \in W$ this Lie algebra is the full
Lie algebra of symmetries of an orthogonal or symplectic form built from 
$\langle -, -\rangle_{\overline{B}_x}$.
\end{remark}

\subsection{Relative theory} We now turn to the relative Hodge theory
of Soergel bimodules. In part this theory is motivated by the unimodality property
\eqref{pos4} of the structure constants $\mu_{x,y}^z$ of
multiplication in the Kazhdan-Lusztig basis.
For all $x, y \in W$ we can find an isomorphism in $\bim$
\begin{equation}
  \label{decomp1}
B_xB_y \cong \bigoplus_{z \in W} V_{x,y}^z \otimes_{\RM} B_z  
\end{equation}
for some graded vector space $V_{x,y}^z$. The structure
constant $\mu_{x,y}^z$ is equal to the graded dimension of the vector
spaces $V_{x,y}^z$. Recall that the unimodality property \eqref{pos4} is equivalent to the
fact that $\mu_{x,y}^z$ is the character of a finite dimensional
$\mathfrak{sl}_2$-module. Of course this would be the case if we can
establish that $V_{x,y}^z$ is actually an
$\mathfrak{sl}_2$-module. This is
equivalent to the existence of an operator $L : V_{x,y}^z \to
V_{x,y}^z$ of degree 2 which satisfies hard Lefschetz (see Remark \ref{rem:sl2}).

The problem is that the decomposition \eqref{decomp1} is not
canonical, and hence it is difficult to produce endomorphisms of the
multiplicity spaces $V_{x,y}^z$. This problem is overcome as follows:
any Soergel bimodule $B$ carries a canonical increasing \emph{perverse
  filtration}  (see \cite[\S 6.2]{EWHodge})
\[
\dots \into \tau_{\le i}(B) \into \tau_{\le i+1}(B) \into \dots
\]
such that all maps are split inclusions of Soergel bimodules (i.e.
the filtration is non-canonically split). The
associated graded (non-canonically isomorphic to $B$) admits a canonical isotypic decomposition
\[
 \gr (B) := \bigoplus \tau_{\le i}(B) / \tau_{\le i-1}(B) =
 \bigoplus_{z \in W}
 H_z(B) \otimes_{\RM} B_z
\]
for certain graded vectors spaces $H_z(B)$.
Moreover any degree $d$
endomorphism of $B$ induces a degree $d$ endomorphism of each
$H_z(B)$.

Applying this construction to the tensor product $B_xB_y$, we may
redefine the vector spaces $V_{x,y}^z$  above as follows:
\[
V_{x,y}^z := H_z(B_xB_y).
\]
These graded vector spaces vanish in degrees of parity different from
$\ell(x) + \ell(y) + \ell(z)$ and are equipped with the following structure:
\begin{enumerate}
\item A graded, symmetric, non-degenerate form induced
  by an ``intersection form'' on $B_xB_y$ (see \cite[\S 2.2]{EWrel}).
\item The structure of a graded $R$-module defined as
  follows: any $r \in R$ of degree $d$ gives a morphism
  \begin{equation}
    \label{eq:3}
    B_xB_y \to B_xB_y(d) : b\otimes b' \mapsto b r \otimes b' = b
    \otimes rb'
  \end{equation}
and hence induces a degree $d$ endomorphism of $V_{x,y}^z$.
\end{enumerate}

The main theorem of \cite{EWrel} is the following:

\begin{theorem} \label{rel}
  For all $x, y, z \in W$ and $\lambda \in \hspos$ the action of $\lambda$ on
  $V_{x,y}^z$ satisfies hard Lefschetz and the Hodge-Riemann bilinear relations.
\end{theorem}

By the discussion above, this theorem immediately implies the
unimodality property \eqref{pos4} of the structure constants of the
Kazhdan-Lusztig basis.

\begin{remark}
  In geometric settings Theorem \ref{rel} can be deduced from the
  relative hard Lefschetz theorem \cite{BBD} and the relative
  Hodge-Riemann bilinear relations \cite{dCM2}.
\end{remark}

\begin{remark}
  Suppose that $W$ is finite with longest element $w_0$. Then one has
  a canonical isomorphism
\[
\overline{B}_x \cong V_{x, w_0}^{w_0}
\]
compatible with forms and the $R$-module structure. (This is a
categorification of the fact that if $\chi$ denotes the ``trivial''
character of $H$, i.e. $H_s \mapsto v^{-1}$ for all $s \in S$,
then we have $b_xb_{w_0} = \chi(b_x)b_{w_0}$ for all $x \in W$.)
In this case Theorem \ref{global} is a special case of Theorem \ref{rel}.
\end{remark}

\begin{remark}
  Theorem \ref{rel} can be used to prove that certain tensor
  categories associated by Lusztig (see \cite{Lust1,Lust2}) to any
  two-sided cell in $W$ are rigid \cite{EWrel}.
\end{remark}

\subsection{Local theory}
We finish our discussion with a brief overview of the local theory,
contained in \cite{WilLoc}. Here the motivation cannot be given
strictly in terms of positivity properties in the Hecke algebra
(although a better understanding of the local route might 
yield an alternative proof of Soergel's conjecture, see \cite[\S 7.5]{WilLoc}).

\begin{remark}
  In geometric settings, the global and relative theory discussed
  above have ``easy'' translations into statements about intersection
  cohomology (see Remarks \ref{rem:hc} and \ref{rem:smod}). In the local setting the
  translation (via the ``fundamental example'' of Bernstein-Lunts
  \cite{BL}) is available, but is more complicated. We will not
  comment on this further below and instead refer the reader to the
  introduction of \cite{WilLoc}, where the connection is discusssed in
  detail.
\end{remark}

For a Soergel bimodule $B$ and $y \in W$, let $\Gamma_y(B)$ (resp. $\Gamma^y(B)$)
denote the largest submodule (resp. quotient) on which one has the
relation $m \cdot r = y(r) \cdot m$ for all $r \in R$. (Thus, for
example, $\Gamma_{\id}(B)$ and $\Gamma^{\id}(B)$ are the Hochschild 
cohomology and homology of the bimodule $B$.)

Both $\Gamma_y(B)$ and $\Gamma^y(B)$ are free graded $R$-modules. The
evident morphism
\[
i_y : \Gamma_y(B) \to \Gamma^y(B)
\]
is injective and becomes an isomorphism if we tensor with $R[1/\alpha_s]_{s \in S}$.

Any $\mu^\vee \in \hg$ determines a specialisation $R \to \RM[z]$
(``restriction to the line $\RM \mu^\vee \subset \hg$'') and if
$\langle \mu^\vee, \alpha_s \rangle \ne 0$ for all $s \in S$ then $\RM[z] \otimes_R i_y$
is an inclusion of graded $\RM[z]$-modules. For any $x \in W$ we
define an $\RM[z]$-module $H_{y,x}^\mu$ via the short exact sequence
\[
0 \to \RM[z] \otimes_R \Gamma_y(B_x) \to \RM[z] \otimes_R
\Gamma^y(B_x) \to H_{y,x}^\mu(1) \to 0.
\]
If we write 
$v^{\ell(x) - \ell(y)}p_{y,x} = \sum_{i \ge 1} a_{y,x}^i v^i$ then the graded rank of
$H_{y,x}^\mu$ is
\[
\sum_{i \ge 1} a_{y,x}^i [i]
\]
(see \eqref{eq:gauss} for the definition of $[i]$).
In particular it vanishes in degrees of the same parity as $\ell(x)$. The
intersection form $B_x$ induces a symmetric non-degenerate form on
$H_{y,x}^\mu$. The main theorem of \cite{WilLoc} is the following:

\begin{theorem} \label{local}
  For any $x < y \in W$ and any $\mu^\vee \in \hpos$, multiplication by $z$
  on $H_{y,x}^z$ satisfies hard Lefschetz and the Hodge-Riemann
  bilinear relations.
\end{theorem}

\begin{remark}
  The main motivation for establishing Theorem \ref{local} is that (by
  work of Soergel \cite{SoeA} and K\"ubel \cite{Ku1,Ku2}) it gives an
  algebraic proof of the Jantzen conjectures on
  the Jantzen filtration on Verma modules for complex semi-simple Lie
  algebras. The first proof of the Jantzen conjectures was given by
  Beilinson and Bernstein \cite{BB}.
\end{remark}

\begin{remark}
The statements and proofs in the local case are less intuitive than
in the previous two settings (global and relative). This might be because we do not yet
have a good framework for discussing the hard Lefschetz theorem and
Hodge-Riemann bilinear relations in equivariant cohomology.
\end{remark}

\section*{Acknowledgements} I would like to thank B. Elias, with whom
it is a pleasure to work on these and related questions. I would also like
to thank A. Beilinson, R. Bezrukavnikov, T. Braden and J. Huh for interesting conversations about
``combinatorial Hodge theory''. Finally, thanks to the referees for
their feedback.

%\section{References}

% It follows a list of references showing you the style in which
% books and journal articles should be listed.

% \renewcommand{\refname}{}    %%%% for this example
% \vspace*{-36pt}              %%%% file only!

\end{document}